\documentclass{amsart}[12pt]
\def\doctype{preprint}

\usepackage{latexsym,amssymb,amsmath}
\usepackage{fancyhdr}
\usepackage{lipsum}
\usepackage{lineno}
\usepackage{hyperref}

\newcommand{\comment}[1]{}

\numberwithin{equation}{section}

\fancypagestyle{titlepage}{
\fancyhead[R]{\doctype}
\fancyhead[CL]{}
\cfoot{\vspace{5pt} \thepage}
}

\newcommand\Section[1]{\section{\bf #1}}

\newtheoremstyle{theorem}
  {12pt}		  
  {0pt}  
  {\sl}  
  {\parindent}     
  {\bf}  
  {. }    
  { }    
  {}     
\theoremstyle{theorem}
\newtheorem{thm}{Theorem}[section]  
\newtheorem{cor}[thm]{Corollary}

\newtheorem{cons}[thm]{Construction}
\newtheoremstyle{definition}
  {12pt}		  
  {0pt}  
  {}  
  {\parindent}     
  {\bf}  
  {. }    
  { }    
  {}     
\theoremstyle{definition}
\newtheorem{ex}[thm]{Example}

\renewcommand\proof{\noindent\textsl{Proof. }}

\makeatletter
\renewcommand*\@maketitle{%
  \normalfont\normalsize
  \@adminfootnotes
  \@mkboth{\@nx\shortauthors}{\@nx\shorttitle}%
  \global\topskip42\p@\relax 
  \@settitle
  \ifx\@empty\authors \else {\vskip 1em
\vtop{\centering\shortauthors\@@par}} \fi
  \ifx\@empty\@date \else {\vskip 1em \vtop{\centering\@date\@@par}}\fi 
  \ifx\@empty\@dedicatory
  \else
    \baselineskip18\p@
    \vtop{\centering{\footnotesize\itshape\@dedicatory\@@par}%
      \global\dimen@i\prevdepth}\prevdepth\dimen@i
  \fi
  \@setabstract
  \normalsize
  \if@titlepage
    \newpage
  \else
    \dimen@34\p@ \advance\dimen@-\baselineskip
    \vskip\dimen@\relax
  \fi
} 
\renewcommand*\@adminfootnotes{%
  \let\@makefnmark\relax  \let\@thefnmark\relax
  \ifx\@empty\@subjclass\else \@footnotetext{\@setsubjclass}\fi
  \ifx\@empty\@keywords\else \@footnotetext{\@setkeywords}\fi
  \ifx\@empty\thankses\else \@footnotetext{%
    \def\par{\let\par\@par}\@setthanks}%
  \fi
\thispagestyle{titlepage}
}
\makeatother

\begin{document}

\title[MOLS]{\large A three-factor product construction for mutually\\ orthogonal latin squares}

\dedicatory{In memory of Rolf Rees.}

\author{Peter J.~Dukes}
\address{\rm Peter J.~ Dukes:
Mathematics and Statistics,
University of Victoria, Victoria, Canada
}
\email{dukes@uvic.ca}

\author{Alan C.H.~Ling}
\address{\rm Alan C.H.~Ling:
Computer Science,
University of Vermont,
Burlington, VT,
USA
}
\email{aling@cems.uvm.edu}

\thanks{Research of Peter Dukes is supported by NSERC}

\date{\today}

\begin{abstract}
It is well known that mutually orthogonal latin squares, or MOLS, admit a (Kronecker) product construction.  We show that, under mild conditions, `triple products' of MOLS can result in a gain of one square.  In terms of transversal designs, the technique is to use a construction of Rolf Rees twice: once to obtain a coarse resolution of the blocks after one product, and next to reorganize classes and resolve the blocks of the second product.  As consequences, we report a few improvements to the MOLS table and obtain a slight strengthening of the famous theorem of MacNeish.
\end{abstract}

\maketitle
\hrule

\Section{Introduction}

A \emph{latin square} is an $n \times n$ array with entries from an $n$-element set of symbols such that every row and column exhausts the symbols (with no repetition).  Often the symbols are taken to be from $\{1,\dots,n\}$.  The integer $n$ is called the \emph{order} of the square.

Two latin squares $L=(L_{ij})$ and $L'=(L'_{ij})$ of order $n$ are \emph{orthogonal} if $\{(L_{ij},M_{ij}): 1 \le i,j \le n\}=\{1,\dots,n\}^2$; that is, two squares are orthogonal if, when superimposed, all ordered pairs of symbols are distinct.

A family of pairwise orthogonal latin squares is normally called \emph{mutually orthogonal latin squares}, and abbreviated `MOLS'.  The maximum size of a family of MOLS of order $n$ is denoted $N(n)$.  It is easy to see that $N(n) \le n-1$ for $n>1$,  with equality if and only if there exists a projective plane of order $n$.  Consequently, $N(q)=q-1$ for prime powers $q$.

Given latin squares $L$, of order $m$, and $M$, of order $n$, their Kronecker product $L \otimes M$ is a latin square of order $mn$.  If $\{L^{(1)},\dots,L^{(k)}\}$ and $\{M^{(1)},\dots,M^{(k)}\}$  are families of MOLS of orders $m$ and $n$, then $\{L^{(1)} \otimes M^{(1)},\dots,L^{(k)} \otimes M^{(k)}\}$ is a family of MOLS of order $mn$.  Hence, $N(mn) \ge \min\{N(m),N(n)\}$.  Combining this with the above remarks yields a `basic' lower bound on $N(n)$.

\begin{thm}[MacNeish's Theorem]
If $n=q_1\dots q_t$ is factored as a product of powers of distinct primes, then $N(n) \ge \min\{q_i-1:i=1,\dots,t\}$. 
\end{thm}

Although it has long been known \cite{CES} that $N(n) \rightarrow \infty$ (in fact, $N(n) \ge n^\frac{1}{14.8}$ for large $n$ is shown in \cite{Beth}), MacNeish's Theorem remains the best known result for many values of $n$, particularly when $n$ has a small number of prime power factors about the same size.  Our main result is a small improvement directed at these challenging cases.

\Section{Transversal Designs and Resolvability}

In what follows, it is convenient to reformulate our discussion of MOLS using the language of designs.

A \emph{transversal design} TD$(k,n)$ consists of an $nk$-element set of \emph{points} partitioned into $k$ \emph{groups}, each of size $n$, and equipped with a family of $n^2$ \emph{blocks} of size $k$ which are pairwise disjoint transversals of the partition.  We have the existence of $r$ MOLS of order $n$ if and only if a TD$(r+2,n)$ exists.  The equivalence is seen by indexing groups of the partition by rows, columns, and symbols from each square. 

In a TD$(k,n)$, a \emph{parallel class} of blocks is a set of $n$ blocks which partition the points. 
If the blocks can be resolved into $n$ parallel classes, such a transversal design is called \emph{resolvable} and denoted RTD$(k,n)$.  The blocks of each parallel class in an RTD$(k,n)$ can be extended by one new point in an extra group.  In this way, it is easy to see that an RTD$(k,n)$ is equivalent to a TD$(k+1,n)$.

More generally, a $\sigma$-\emph{parallel class} is a configuration of  blocks which covers every point exactly $\sigma$ times.  For a list of positive integers $\Sigma=[\sigma_1,\dots,\sigma_t]$ summing to $n$, a TD$(k,n)$ is $\Sigma$-resolvable if the blocks can be resolved into $\sigma_i$-parallel classes for $i=1,\dots,t$. In writing a list $\Sigma$, we use `exponential notation' such as $\sigma^j$ to abbreviate $j$ occurrences of $\sigma$.

Let us say that a TD$(k,m)$ admits a $(\sigma,\gamma)$-\emph{group partition} if each of the groups of size $m$ is written on some group $G$, if there exists a subset $H$ of $G$ with $|H|=\sigma$ and there exists a partition $\mathfrak{b}$ of the blocks so that, for every class $\mathcal{B} \in \mathfrak{b}$, the set $\{H*B:B \in \mathcal{B} \}$ is a $\gamma$-parallel class.  We should interpret $H*B$ as $|H|$ blocks obtained under action of $H$ on $B$.

A TD$(k,m)$ always admits two `trivial' $(\sigma,\gamma)$-group partitions at each of two extremes.  We have a $(1,m)$-group partition, in which $H$ is the identity subgroup and $\mathfrak{b}$ is the trivial partition with all blocks in the same class, and also an $(m,1)$-group partition, in which $H=G$ is the full group and $\mathfrak{b}$ is the partition into singleton block classes.  An RTD$(k,m)$ admits a $(1,1)$-group partition.

The following is a special case of Construction 2 in \cite{Rees-prod}, due to Rolf Rees.  A similar construction by the same author later appears in \cite{Rees-ttd} in the context of transversal designs.

\begin{cons}[Rees,\cite{Rees-ttd,Rees-prod}]
\label{Rees}
Let $\Sigma=[\sigma_1,\dots,\sigma_t]$.
Suppose there exists a $\Sigma$-resolvable TD$(k,n)$ and a  TD$(k,m)$ admitting, for each $i$, a $(\sigma_i,\gamma_i)$-group partition.
Then there exists a $\Gamma$-resolvable TD$(k,mn)$, where $\Gamma$ consists of $m\sigma_i/\gamma_i$ copies of $\gamma_i$, for $i=1,\dots,t$.
\end{cons}

This is essentially the standard direct product in which blocks of the TD$(k,n)$ are replaced by copies of the TD$(k,m)$.  For resolving, the key idea is this: given a $\sigma$-parallel class $\mathcal{C}$ in the TD$(k,n)$ and a $(\sigma,\gamma)$-group partition $(H,\mathfrak{b})$ of the TD$(k,m)$, we can `split' the occurrences of blocks in $\mathcal{C}$ incident with each point $x$ using bijections onto $\{x\} \times H$ in the product.  Then, we can select $\gamma$-parallel classes in the product according to the action of $H$ on $\mathfrak{b}$.  

Note that a $(\sigma,\gamma)$-group partition is a stronger hypothesis than is actually needed for the construction in \cite{Rees-prod}, where different subsets of $G$ can be taken in each group.  Since we only need the `easy' group partitions mentioned above, we adopt our stronger, simplified hypothesis for clarity.  In any case, with this construction we are ready to state and prove our main result. 

\begin{thm}
\label{triple}
For integers $a,b,c$ with $a \le b \le c$, we have
$$N(abc) \ge \min\{N(a)+1,N(b),N(c)\}.$$
\end{thm}

\proof
Equivalently, we construct an RTD$(k,abc)$ given TD$(k,a)$, RTD$(k,b)$, and RTD$(k,c)$.  Put $c=aq+r$, $0 \le r < a$.
The last ingredient, an RTD$(k,c)$, can be regarded instead as a $\Sigma$-resolvable TD$(k,c)$, where $\Sigma=[1^r,a^q]$. Apply Construction~\ref{Rees}, using $(1,a)$- and $(a,1)$-group partitions of a TD$(k,a)$.  The result is a $\Gamma$-resolvable TD$(k,ac)$, where
$$\Gamma=\left[1^{qa^2},a^r \right].$$
Since $qa^2+ar = ac > rb$, we may reorganize the classes to obtain a $\Gamma'$-resolvable TD$(k,ac)$, where
$$\Gamma'=\left[1^{qa^2-r(b-a)},b^r \right]= \left[1^{ac-rb}, b^r \right].$$
Now, apply Construction~\ref{Rees} again with an RTD$(k,b)$, using $(1,1)$- and $(b,1)$-group partitions.  
\qed

\begin{cor}
\label{pqr}
For prime powers $p \le q \le r$, we have $N(pqr) \ge p.$
\end{cor}

In particular, we have `easy' proofs of $N(18),N(30) \ge 2$.  In fact, Rees had already obtained those orthogonal latin squares using his construction, bypassing the step of `reorganizing' classes.  This step is the key contribution of this paper; we next provide an example illustrating its usefulness.

\begin{ex}
Consider $p=8$, $q=9$, $r=13$.  There exists an RTD$(9,13)$.  By amalgamating parallel classes, this can be viewed also as a $[1^5,8]$-resolvable TD$(9,13)$.
Now, consider a TD$(8,8)$, which admits both an $(8,1)$- and a $(1,8)$-partition.  It follows by Construction~\ref{Rees} that there exists a $[1^{64},8^5]$-resolvable TD$(9,104)$.  Again, by restructuring classes, this can be viewed instead as $[1^{59},9^5]$-resolvable.  Since there exists an RTD$(9,9)$, it admits both $(1,1)$- and $(9,1)$-group partitions.  Applying Construction~\ref{Rees} again yields a $[1^{9 \times 59},1^{5 \times 81}]$-resolvable TD$(9,936)$.  In other words, we have an RTD$(9,936)$ or equivalently a TD$(10,936)$.  This gives $N(936) \ge 8$.
\end{ex}

\begin{table}[h!]
$$\begin{array}{rcr|c|c}
\text{factorization} & & n~ & N_\text{HCD}(n) & N(n) \ge\\
\hline
8 \times 9 \times 13 &=& 936 & 7 & 8 \\
8 \times 9 \times 17 &=& 1224 & 7 & 8 \\
8 \times 11 \times 13&=&  1144 & 7 & 8\\
16 \times 17 \times 19 &=& 5168 & 15 & 16\\
16 \times 17 \times 25 &=& 6800 & 15 & 16\\
16 \times 19 \times 31 &=& 9424 & 15 & 16\\
17 \times 19 \times 23 &=& 7429 & 16 & 17\\
17 \times 19 \times 29 &=& 9367 & 16 & 17 
\end{array}$$
\caption{Improvements to the table of MOLS in \cite{Handbook} from Theorem~\ref{triple}.}
\label{mols-improvements}
\end{table}

Table~\ref{mols-improvements} gives a list of similar lower bounds on $N(n)$ which improve upon the bounds $N(n) \ge N_\text{HCD}(n)$ stated in \cite{Handbook} for $n<10,000$.  In fact, the entries of our table account for all integers $n$ having three or more prime power exact divisors such that $N_\text{HCD}(n)$ reports the MacNeish bound.

We should remark that $p,q,r$ need not be assumed relatively prime in Corollary~\ref{pqr}.  For instance, $N(31 \times 2^{t}) \ge 31$ for $t \ge 10$ by factoring $2^{t}=2^{5} \times 2^{t-5}$.  However, we could find no cases for $n<10,000$ where splitting prime powers improves the state of the art.  We fail to improve $N_\text{HCD}(17 \cdot 2^9)= 16$, for instance, since $2^4<17$.

In closing, we should also note that Theorem~\ref{triple} can be iterated, though in light of $N(n) \ge n^{\frac{1}{14.8}}$ for large $n$, it is not (at least asymptotically) worthwhile to iterate very often.  But, for example, since $N(8 \times 9 \times 13) \ge 8$ and $8 \times 9 \times 13 < 5^5$, we also have $N(8 \times 9 \times 13 \times 5^{10}) \ge 9$.

\end{document}